\documentclass[11pt]{amsart}
\pdfoutput=1
\usepackage{bm}
\usepackage{graphicx}
\usepackage{mathpazo}
\usepackage[T1]{fontenc}
\usepackage{latexsym}   
\usepackage{amssymb}    
\usepackage{mathrsfs}   
\usepackage{amsmath}    
\usepackage{amsthm}     
\usepackage[margin=1in]{geometry}
\usepackage[color,cmtip,all,matrix,arrow,tips,curve]{xy}

\usepackage[colorlinks,breaklinks=true]{hyperref}

\usepackage{multibib} 
\newcites{A}{Works that cite Manin (1961) or Yui (1978)} 

\numberwithin{equation}{section}

\let\stdthebibliography\thebibliography
\let\stdendthebibliography\endthebibliography
\renewenvironment*{thebibliography}[1]{%
  \stdthebibliography{100}}
  {\stdendthebibliography}

\usepackage{totcount}
\newtotcounter{appbibcounter}

\theoremstyle{definition}

\newcommand{\half}[1]{\frac{#1}{2}}
\newcommand{\st}[1]{\left\{#1\right\}}

\DeclareMathOperator{\charpoly}{charpoly}
\DeclareMathOperator{\id}{id}
\DeclareMathOperator{\Jac}{Jac}
\DeclareMathOperator{\mat}{Mat}
\DeclareMathOperator{\residue}{res}

\newcommand{\derham}{{\mathrm{dR}}}
\newcommand{\cartier}{{\mathcal C}}
\newcommand{\frobenius}{{\mathcal F}}

\newcommand{\calB}{{\mathcal B}}
\newcommand{\calC}{{\mathcal C}}
\newcommand{\calD}{{\mathcal D}}
\newcommand{\calF}{{\mathcal F}}
\newcommand{\calO}{{\mathcal O}}
\newcommand{\calR}{{\mathcal R}}

\newcommand{\ff}{{\mathbb F}}
\newcommand{\Z}{{\mathbb Z}}
\newcommand{\fbar}{\bar{\ff}}

\newcommand{\inv}{^{-1}}
\newcommand{\cross}{\times}
\newcommand{\iso}{\cong}

\newcommand{\mybar}[1]{
  \mathchoice
  {#1\llap{$\overline{\phantom{\displaystyle\rm#1}}$}}
  {#1\llap{$\overline{\phantom{\textstyle\rm#1}}$}}
  {#1\llap{$\overline{\phantom{\scriptstyle\rm#1}}$}}
  {#1\llap{$\overline{\phantom{\scriptscriptstyle\rm#1}}$}}
}  
\renewcommand{\bar}{\mybar}

\title[Hasse--Witt and Cartier--Manin matrices]{Hasse--Witt and Cartier--Manin matrices: \\ A warning and a request}
\author{Jeffrey D. Achter}
\address{Department of Mathematics, Colorado State University, Fort
Collins, CO 80523} 
\urladdr{\url{http://www.math.colostate.edu/~achter}}
\email{\href{mailto:j.achter@colostate.edu}{j.achter@colostate.edu}}
\thanks{JDA's work partially supported by NSA grant H98230-16-1-0046}
\author{Everett W. Howe}
\address{Center for Communications Research, 4320 Westerra Court, San Diego, CA 92121-1967}
\curraddr{Unaffiliated, San Diego}
\urladdr{\url{http://ewhowe.com}}
\email{\href{mailto:however@alumni.caltech.edu}{however@alumni.caltech.edu}}
\date{5 February 2020}

\begin{document}

\begin{abstract}
Let $X$ be a curve in positive characteristic.
A \emph{Hasse--Witt matrix} for $X$ is a matrix that represents
the action of the Frobenius operator on the cohomology group $H^1(X,\calO_X)$
with respect to some basis.
A \emph{Cartier--Manin matrix} for $X$ is a matrix that represents
the action of the Cartier operator on the space of holomorphic differentials of~$X$
with respect to some basis.
The operators that these matrices represent are adjoint to one another, so 
Hasse--Witt matrices and the Cartier--Manin matrices are related to one another, 
but there seems to be a fair amount of confusion in the literature about
the exact nature of this relationship.
This confusion arises from differences in terminology, 
from differing conventions about whether matrices act on the left or on the right,
and from misunderstandings about the proper formul\ae{} for 
iterating semilinear operators.
Unfortunately, this confusion has led to
the publication of incorrect results. In this paper we present the issues
involved as clearly as we can, and we look through the literature to see
where there may be problems. We encourage future authors to clearly
distinguish between Hasse--Witt and Cartier--Manin matrices, in the
hope that further errors can be avoided.

\end{abstract}

\maketitle

\section*{Prologue}
\subsection*{An example}

Consider the genus-$2$ hyperelliptic curve $X$ over $\ff_{125}$ with affine
model
\begin{equation}
\label{eq:example}
y^2 = f(x) = x^5+x^4+\alpha^{92}x^3+\alpha^{18}x^2+\alpha^{56}x\,,
\end{equation}
where $\alpha \in \ff_{125}$ satisfies $\alpha^3+3\alpha+3 = 0$.  
Let us compute the $5$-rank of (the Jacobian of) $X$.

On one hand, we can follow Yui~\cite{Yui1978} and compute the effect of the Cartier operator on the space of regular one-forms.  
Let $c_m$ be the coefficient of $x^m$ in $f(x)^{(5-1)/2}$. Yui~\cite[p.~381]{Yui1978} constructs a matrix
(denoted $A$ in her paper, but denoted here by $Y$ to prevent a conflict of notation later on) given by
\begin{align*} 
Y 
&= \begin{pmatrix}
      c_{5\cdot 1 - 1} & c_{5\cdot 1-2} \\
      c_{5\cdot 2 - 1} & c_{5\cdot 2-2}
   \end{pmatrix} \\
&= \begin{pmatrix}
      \alpha^{41} & \alpha^{105} \\
      2           & \alpha^{95}
   \end{pmatrix}.
\intertext{We compute as well that the image of $Y$ under 
the Frobenius automorphism $\sigma$ of~$\ff_{125}$ is given by}
Y^\sigma 
&= \begin{pmatrix}
      c_{5\cdot 1 - 1}^\sigma & c_{5\cdot 1-2}^\sigma \\
      c_{5\cdot 2 - 1}^\sigma & c_{5\cdot 2-2}^\sigma
   \end{pmatrix}\\
&= \begin{pmatrix}
      \alpha^{81} & \alpha^{29} \\
      2 & \alpha^{103}
   \end{pmatrix},
\intertext{and the product $Y \cdot Y^\sigma$ is}
Y \cdot Y^\sigma 
&= \begin{pmatrix}
      \alpha^{32} & \alpha^{104} \\
      \alpha^{22} & \alpha^{94}
   \end{pmatrix}.
\end{align*}
Since this last matrix has rank one, according to 
Yui's Lemma E~\cite[p.~387]{Yui1978}
we should be able to conclude that $X$ has $5$-rank one.

On the other hand, $X$ is actually supersingular; indeed, its L-polynomial is 
$(1+125T^2)^2$, and thus the only slope of its normalized $5$-adic Newton
polygon is $1/2$.  In particular, $X$ has $5$-rank zero.

Our aim in this note is to tease out the source of this dissonance.

\subsection*{Genesis of this project}

We noticed this discrepancy while attempting to obtain numerical data
in support of some earlier work~\cite{AchterHowe2017}.  Moreover, we
found that one of us invoked an erroneous
formula in a separate project~\citeA{Howe2008} (see Section~\ref{ss:genus2}).

Works such as Yui's 1978 paper~\cite{Yui1978}, as well as its antecedents 
(including works by Manin~\cite{Manin1961,Manin1961english}) and consequents,
rely on the construction and analysis
of certain semilinear operators.  Since the work of 
Hasse and Witt~\cite{HasseWitt1936},
it has been understood that there is such an operator,
acting on some subquotient of the de~Rham cohomology of a given curve
$X$ in characteristic~$p$, that encodes information about the
$p$-torsion group scheme of the Jacobian of $X$.  The ideas of Hasse
and Witt are beautifully clear, but one must navigate around several
potential sources of error in order to arrive at a correct formula.
Indeed, one must decide whether to work with the summand
$H^0(X,\Omega^1_X)$ or the quotient $H^1(X,\calO_X)$ of
$H^1_\derham(X)$; this choice, in turn, determines whether the
operator in question is $\sigma$-linear or
$\sigma\inv$-linear, where $\sigma$ is the $p$-th powering map on
the base field.  One is given a further opportunity to make a
``sign error'' when one chooses bases for these vector spaces
and then decides whether the semilinear operator acts on the right
or on the left.\footnote{Of course, there is no literal ``sign'' to get wrong 
          in any of the formul\ae{} we discuss, but the terminology is
          suggestive of the fact that two such errors will typically cancel
          one another out. We will continue 
          to use the the term ``sign error'' in this sense throughout the paper.
}
Given these multiple opportunities for mistake, it is hardly surprising
that there are occasional misstatements in the literature.

Conversations with others suggest to us that the community has an 
interest in (re)documenting these semilinear methods, especially in view of
the continuing expansion of the role of computing in arithmetic geometry.
With this backdrop, we offer the following survey
of Cartier--Manin and Hasse--Witt matrices.

In Section~\ref{sec:matrices} we review basic facts about the representation
of semilinear operators by matrices.
In Section~\ref{sec:HW&CM} we define the Cartier operator on the space of
holomorphic differentials of a curve $X$ and the Frobenius operator on the
cohomology group $H^1(X,\calO_X)$, in its guise as a quotient group of the space
of r\'epartitions of~$X$. 
The Cartier--Manin and Hasse--Witt matrices represent these two operators. 
In Section~\ref{sec:formulaire} we follow the work of
Manin~\cite{Manin1962,Manin1962english} and Yui~\cite{Yui1978} to explicitly
calculate the Cartier--Manin matrix of a hyperelliptic curve, and we 
resolve the problem posed by the example in our Prologue.
In Section~\ref{sec:Manin&Yui} we review the papers of Manin and Yui, 
keeping a watchful eye out for sign errors. 
We close in Section~\ref{sec:subsequent} with a review of the 
literature that cites Manin and Yui, to see whether any
sign errors have propagated. Fortunately, there are only a few papers that 
contain results or examples that are in error.

Of course, it is unpleasant to find any errors at all in published papers.
We have a suggestion for the community, which we hope will help prevent 
this type of sign error in the future: Please be careful with terminology.
If you are working with the Cartier operator on differentials, refer to the
matrix representation as the \emph{Cartier--Manin} matrix; if you are working
with the Frobenius operator on $H^1(X,\calO_X)$, refer to the
matrix representation as the \emph{Hasse--Witt} matrix.  These matrices 
are \emph{related} to one another, but they are not \emph{equal} to one another,
and they represent semilinear operators with different properties.

\subsection*{Acknowledgments}
We thank 
Yuri Manin, 
Noriko Yui, 
Arsen Elkin,
Pierrick Gaudry,
Takehiro Hasegawa, 
Rachel Pries,
Andrew Sutherland,
Saeed Tafazolian,
Doug Ulmer,
Felipe Voloch,
and Yuri Zarhin,
as well as the referees,
for their comments on draft versions of this paper.
We thank Dino Lorenzini and Felipe Voloch for drawing our attention to a mistake
in Section~\ref{subseczeta} that appeared in the published version of this paper
(corrected below); the irony does not escape us.

\section{Matrices and semilinear algebra}
\label{sec:matrices}

We start with some notation concerning the use of matrices
to represent semilinear algebra.

Let $K$ be a field; we work exclusively with finite-dimensional
$K$-vector spaces. 

\subsection{Bases, matrices, and linear operators}

Let $W$ be a vector space with basis $\calC = \st{w_1, \cdots, w_n}$.
Any $w\in W$ is expressible as $w = \sum c_i w_i$; let $[w]_{\calC}$
denote the corresponding column vector  
\[
[w]_{\calC} = \begin{pmatrix}c_1\\c_2\\\vdots\\c_n\end{pmatrix}.
\]
Now let $V$ be an $m$-dimensional vector space with chosen basis
$\calB$, and let $f\colon W \to V$ be a linear transformation.  Define
numbers $a_{ij}$ by   
\[
f(w_j) = \sum_{i=1}^m a_{ij}v_i\,.
\]
The matrix of $f$  relative to the chosen bases $\calC$ and $\calB$ is
\[
[f]_{\calB\leftarrow \calC} = (a_{ij}) \in \mat_{m\times n}(K)\,.
\]
Matrix multiplication is defined so that, with this notation,
\[
\left[f(w)\right]_{\calB} = [f]_{\calB \leftarrow \calC}\cdot \left[w\right]_{\calC}\,.
\]

Change of basis works as follows.  Let $f\colon V \rightarrow V$ be an
endomorphism, and let $\calB$ and $\calD$ be two different bases for
$V$.  Then 
\begin{align*} \left[f\right]_{\calD\leftarrow \calD} &= \left[\id\right]_{\calD \leftarrow \calB}
\left[f\right]_{\calB \leftarrow \calB} \left[\id\right]_{\calB \leftarrow \calD}\,; 
\intertext{if $S =[\id]_{\calB \leftarrow \calD}$, then} 
\left[f\right]_{\calD\leftarrow \calD}&=S\inv \left[f\right]_{\calB \leftarrow \calB} S\,.
\end{align*}
(Of course, if one prefers that matrices act on the right, then one
consistently writes elements of the vector space as row vectors, and
the matrix that represents the action of a linear operator is the
\emph{transpose} of the matrix described above.)

\subsection{Semilinear algebra}
\label{subsecsemilinear}

Let $\epsilon$ be an automorphism of $K$.  Now suppose that $f\colon V \to   
V$ is $\epsilon$-linear, in the sense that for $a \in K$ and $v\in V$,
\[
f(av) = a^\epsilon f(v)\,.
\]

Naturally, $f$ is determined by its effect on a basis, but the use of
the matrices changes a little bit.  Let $\calB = \st{v_1, \cdots,
v_n}$ be a basis, and again define numbers $a_{ij}$ by
\[
f(v_j) = \sum_i a_{ij}v_i\,.
\]
If $v = \sum_j c_j v_j$ then
\[
f(v) = \sum_j f(c_jv_j) 
     = \sum_j c_j^\epsilon f(v_j) 
     = \sum_j \Big(\sum_i a_{ij}v_i\Big) c_j^\epsilon\]
and so
\[\left[f(v)\right]_{\calB}
= \left[f\right]_{\calB \leftarrow \calB}\cdot \left[v\right]_{\calB}^\epsilon\,,
\]
where $B^\epsilon$ is the matrix obtained by applying $\epsilon$ to each entry of $B$.

Similarly, change of basis is accomplished with $\epsilon$-twisted conjugacy:
\begin{align*}
\left[f\right]_{\calD\leftarrow \calD} &=
\left[\id\right]_{\calD \leftarrow \calB} \left[f\right]_{\calB \leftarrow \calB} \left[\id\right]_{\calB\leftarrow\calD}^\epsilon \\
&= S\inv \left[f\right]_{\calB \leftarrow \calB} S^\epsilon.
\end{align*}

If we suppress our subscripts for a moment,
then the iterates of $f$ are represented by
\begin{align*}
\left[f\circ f\right] &= \left[f\right] \left[f\right]^\epsilon \\
\left[f^{\circ r}\right] &= \left[f\right] \left[f\right]^{\epsilon} \left[f\right]^{\epsilon^2} \cdots \left[f\right]^{\epsilon^{r-1}}.
\end{align*}

\noindent (Again, if one wants matrices to act on the right, then the highest
iterate of $\epsilon$ is applied to the \emph{leftmost} factor in the
$r$-fold product.)
 
\subsection{Adjointness}
\label{subsecdual}

Let $V^*$ be the dual vector space of $V$ and let $(\cdot,\cdot)\colon V \times V^* \to K$ be the natural pairing.  Continue to let $f\colon V \to V$ be $\epsilon$-linear, and let $\delta=\epsilon\inv$.  
The adjoint $f^*$ of $f$ with respect to the pairing $(\cdot,\cdot)$ is $\delta$-linear and is characterized by the relation
\[
(v,f^*w^*) = (fv,w^*)^{\delta}
\]
for all $v\in V$ and $w^*\in V^*$.
Let $\calB^*=\st{v_1^*, \cdots, v_n^*}$ be the basis dual to $\calB$.  Since for $1 \le j, \ell \le n$
we have
\[
(v_j,f^*v_\ell^*) = (fv_j,v_\ell^*)^{\delta} = \big(\sum_{i} a_{ij}v_i,v_\ell^*\big)^{\delta} = a_{\ell j}^{\delta}\,,
\]
we find that
\[
f^* v_\ell^* = \sum_{j} a_{\ell j}^\delta v_j^*
\]
and therefore
\begin{equation}
\label{eq:dual}
[f^*]_{\calB^*\leftarrow \calB^*} = \left([f]_{\calB\leftarrow\calB}^{\delta}\right)^\intercal\,
\end{equation}
where $^\intercal$ indicates the transpose of a matrix.

\section{Hasse--Witt and Cartier--Manin matrices}
\label{sec:HW&CM}

We record here some properties of the Frobenius and Cartier operators
and their representations by Hasse--Witt and Cartier--Manin matrices,
deferring a complete exposition to, for example, Serre~\cite{Serre1958}. 
Let $k$ be a perfect field of characteristic $p>0$.
Let $\sigma\colon k \to k$ be the Frobenius automorphism, and let $\tau$ be
its inverse.
Finally, let $X/k$ be a smooth, projective curve of genus $g>0$.

\subsection{Cohomology groups}

The Hodge to de~Rham spectral sequence gives a canonical exact
sequence 
\[
\xymatrix{
0 \ar[r] &  H^0(X,\Omega^1_X) \ar[r] & H^1_\derham(X)  \ar[r]&
H^1(X,\calO_X) \ar[r] & 0\,.
}
\]
There is a canonical duality between the $g$-dimensional $k$-vector
spaces $H^0(X,\Omega^1_X)$ and $H^1(X,\calO_X)$.  This duality is
realized by cup product and the trace map:
\[
\xymatrix{H^0(X,\Omega^1_X)\cross H^1(X,\calO_X) \ar[r] & H^1(X,\Omega^1_X) \ar[r]^-\sim & k\,.
}
\]

If $k$ is algebraically closed, Serre~\cite[\S~8]{Serre1958} gives the following explicit
description of this pairing.  Let $\calR = \calR(X)$ be the ring of
r\'epartitions on $X$ --- that is, the subring of $\prod_{P\in X(k)} k(X)$ consisting of those elements $\st{r_P}$ for which, for
all but finitely many $P$, the function $r_P$ is regular at $P$.  Let $\calR(0)$ be
the subring consisting of those r\'epartitions such that each $r_P$ is
regular at $P$.  Then there is an isomorphism
\[
H^1(X,\calO_X) \iso \frac{\calR}{\calR(0)+k(X)}\,,
\]
where we view $k(x)$ as a subring of $\calR$ via the diagonal embedding.
The duality between this space and $H^0(X,\Omega^1_X)$ then admits the
description 
\begin{equation}
\label{eqpairing}
\xymatrix@R-2pc{
H^0(X,\Omega^1_X)\cross H^1(X,\calO_X) \ar[r]& \rlap{$k$}\phantom{\sum_{P\in X(k)} \residue_P(r_P\omega)\,,}\\
{}\phantom{H^0(X,\Omega^1_X)\cross H^1(X,\calO_X)}\llap{$(\omega,r)$} \ar@{|->}[r] & \sum_{P\in X(k)} \residue_P(r_P\omega)\,,
}\end{equation}
where $\residue_P$ denotes the residue at the point~$P$.

\subsection{The Cartier operator and the Cartier--Manin matrix}
\label{subseccartier}
Cartier~\cite{Cartier1957} (see also Katz~\cite[\S 7]{Katz1970})
defines an operator on the de Rham complex of a smooth proper variety of arbitrary
dimension.  In the special case of a curve $X$, this give rise to a map from $H^0(X,\Omega^1_X)$ to itself.  We follow here the explicit description given by Serre~\cite[\S 10]{Serre1958}.

Let $P$ be a closed point on  $X$ and let $t$ be a uniformizing
parameter at $P$.  Then the functions $1, t, \cdots, t^{p-1}$ form a
$p$-basis for the local ring $\calO_{X,P}$, that is, a basis for $\calO_{X,P}$ as a module over $\calO_{X,P}^p$.  Any $1$-form
holomorphic at $P$ admits an expression
\begin{align*}
\omega &= \left(\sum_{j=0}^{p-1} f_j^p t^j\right)\, dt
\intertext{for certain $f_j \in \calO_{X,P}$, and  one declares that}
\cartier(\omega) &= f_{p-1}\, dt\,.
\end{align*}
The value of $\cartier(\omega)$ does not depend on the choice of 
uniformizer~$t$, and the map $\cartier$ can be extended to give
a map $\Omega^1_{k(X)/k} \to \Omega^1_{k(X)/k}\,.$

It is not hard to see that, for $\omega$, $\omega_1$, and $\omega_2$ in
$\Omega^1_{k(X)/k}$ and for $f \in k(X)$, one has
\begin{align*}
\cartier(\omega_1+\omega_2) &= \cartier(\omega_1)+\cartier(\omega_2)\\
\cartier(f^p\omega) &= f\cartier(\omega)\,.
\end{align*}
In particular, the Cartier operator restricts to give a $\tau$-linear
operator 
\[
\xymatrix{
H^0(X,\Omega^1_X) \ar[r]^{\cartier} & H^0(X,\Omega^1_X)\,.
}\]
(Yui~\cite{Yui1978} refers to this as the \emph{modified} Cartier operator.)
A matrix associated to $\cartier$ and a choice of basis for
$H^0(X,\Omega^1_X)$ is called a \emph{Cartier}, or
\emph{Cartier--Manin}, matrix for~$X$.

\subsection{The Frobenius operator and the Hasse--Witt matrix}

There is also a Frobenius operator
\[\xymatrix{
H^1(X,\calO_X) \ar[r]^\frobenius & H^1(X,\calO_X)
}\]
which, under the isomorphism $H^1(X,\calO_X) \iso \calR/(\calR(0)+k(X))$,
takes the class of a  r\'epartition $r = \st{r_P}$ to the class
of $\big\{r_P^p\big\}$. 
In particular, $\frobenius$ is a $\sigma$-linear operator.
Following Serre, we call any matrix associated to $\frobenius$ and a choice
of basis  a \emph{Hasse--Witt}
matrix for $X$.

Like the Cartier operator, the Frobenius operator admits a generalization to varieties of arbitrary dimension.  For a smooth variety for which the Hodge to de Rham spectral sequence degenerates at $E_1$, Katz defines \citeA[(2.3.4.1.3), p.~27]{Katz1972} a $\sigma$-linear operator on each cohomology group of the structure sheaf.
In the special case of a smooth projective hypersurface $Y/k$ of dimension $n$, Katz gives an explicit formula for the action of this operator on $H^n(Y,\mathcal O_Y)$ in terms of a defining polynomial for $Y$ \citeA[Algorithm (2.3.7.14), p.~35]{Katz1972}.

\subsection{Adjointness}

Serre goes on to show \cite[Proposition~9, p.~40]{Serre1958} that $\frobenius$ and $\cartier$ are adjoint with respect to the pairing $(\cdot,\cdot)$ of \eqref{eqpairing}, in the sense (see Section \ref{subsecdual}) that
\begin{equation}
\label{eq:frobcartierdual}
(\omega, \frobenius r) = (\cartier \omega, r)^\sigma.
\end{equation}
By \eqref{eq:dual}, if $B$ is a Cartier--Manin matrix for $X$, then $(B^\sigma)^\intercal$ is a Hasse--Witt matrix for~$X$.  Conversely, if $A$ is a Hasse--Witt matrix for $X$, then $(A^\tau)^\intercal$ is a Cartier--Manin matrix for~$X$.

\subsection{Zeta functions}
\label{subseczeta}
Now suppose $X$ is a curve of genus $g$ over $\ff_{p^e}$, 
the field with $q = p^e$ elements. 
The zeta function of $X$ has the form
\[
Z_{X/\ff_q}(T) = \frac{L(T)}{(1-T)(1-qT)}\,,
\]
where $L(T) \in \Z[T]$ is a polynomial of degree~$2g$. 
The $e$-fold iterate of $\frobenius$ is $\ff_q$-linear, and if $A$ is the
Hasse--Witt matrix for $X$ corresponding to a particular choice of basis,
then $\frobenius^e$ is represented on that same basis 
by $AA^\sigma\cdots A^{\sigma^{e-1}}$.
By \cite[Th\'eor\`eme~3.1]{Katz1973}, we then have 
\begin{equation}
\label{EQ:zeta}
\det(\id-AA^\sigma\cdots A^{\sigma^{e-1}}\,T) \equiv L(T) \bmod p.
\end{equation}
Using \eqref{eq:dual}, we find that $\cartier^{e}$ and $\frobenius^e$ are adjoint
$\ff_q$-linear operators.  In particular, if $B$ is any Cartier--Manin matrix for $X$,
then 
\[
\det(\id-BB^\tau\cdots B^{\tau^{e-1}}\,T) \equiv L(T) \bmod p.
\]
Similarly, the characteristic polynomial $\chi_{X/\ff_q}(T)$ of the relative 
Frobenius endomorphism of $\Jac X$ satisfies
\[
\chi_{X/\ff_q}(T) \equiv T^g \det(T \cdot \id - [\frobenius^{e}])
\equiv T^g \det(T \cdot\id - [\cartier^{e}]) \bmod p\]
(see \cite[Theorem~1]{Manin1961} and \cite[Theorem~1]{Manin1961english}).

Note: Earlier versions of this paper, including the published version,
claimed incorrectly 
that $\charpoly_{\frobenius^e}(T) \equiv L(T) \bmod p.$ The relation we intended
to write is \[{T^g \charpoly_{\frobenius^e}(1/T) \equiv L(T) \bmod p},\] which follows
immediately from~\eqref{EQ:zeta}.

\section{Cartier--Manin matrices for hyperelliptic curves}
\label{sec:formulaire}

We use the methods of Manin~\cite{Manin1962,Manin1962english} and Yui~\cite{Yui1978} 
to give a formula for a Cartier--Manin matrix of a hyperelliptic curve.
We then use this formula to compute such a matrix for the curve \eqref{eq:example}, and independently compute a Hasse--Witt matrix to verify our work.

\subsection{An explicit formula}
\label{subsec:formulaire}
Let $k$ be a perfect field of odd characteristic $p$, and let $X/k$ be a
hyperelliptic curve of genus $g$ with affine equation
$y^2=f(x)$,
where $f(x) \in k[x]$ is square-free of degree $2g+1$ or $2g+2$.

As a basis for $H^0(X,\Omega^1_X)$ we choose
\begin{equation}
\label{eqbasis}
\calB = \st{ \omega_i = x^{i-1}\frac{dx}{y} : 1 \le i \le g}.
\end{equation}
If we write $f(x)^{\half{p-1}} = \sum c_m x^m$,
we obtain the following equalities of differentials on $X$:
\begin{align*}
\frac{dx}{y} &= \frac{(y^2)^{\half{p-1}}}{y^p}dx \\
&= \frac{f(x)^{\half{p-1}}}{y^p}dx\\
&= y^{-p}\bigg( \sum_{m\ge0} c_m x^m \bigg)dx\,.
\end{align*}
We find that
\[
\omega_j = x^{j-1}\frac{dx}{y} = y^{-p}\bigg( \sum_{m\ge0} c_{m}x^{m+j-1} \bigg)dx\,.
\]
If we apply the Cartier operator to $\omega_j$, the only
terms that will make a contribution are the 
terms where $m+j-1\equiv p-1\bmod p$.
In particular, we only need consider $m$ of the form $ip-j$, for $i=1,\ldots,g$. We find that
\begin{align*}
\cartier(\omega_j) 
&= \cartier\bigg(y^{-p} \Big(\sum_{i = 1}^{g}c_{ip-j}x^{ip-p}\Big) x^{p-1}\, dx\bigg) \\
&= \sum_{i = 1}^{g} \cartier\bigg(\left(c_{ip-j}^\tau x^{i-1}/y\right)^p x^{p-1}\, dx\bigg)\\
&= \sum_{i = 1}^{g} c^\tau_{ip-j}x^{i-1}/y\, dx \\
&= \sum_{i\ge 1} c^\tau_{ip-j} \omega_i\,.
\end{align*}
If we let $B \in \mat_g(k)$ be the matrix with entries 
$B_{ij} = c^\tau_{ip-j}$,
then left-multiplication by $B$ calculates the effect of $\cartier$ in the basis $\calB$.

\subsection{The example, revisited}

We reconsider the curve \eqref{eq:example} and the associated matrix $Y$.  Then
\[
B = Y^\tau = \begin{pmatrix}
\alpha^{33} & \alpha^{21}\\
2 & \alpha^{19}
\end{pmatrix}.
\]
We compute the effect of the second iterate of the Cartier operator as 
\[
\null[\cartier^{\circ 2}]_{\calB \leftarrow \calB} = [\cartier] [\cartier]^\tau = B B^\tau =
\begin{pmatrix}
0&0\\0&0
\end{pmatrix};
\]
this reflects the supersingularity of our original curve.

For thoroughness, we will use direct computation to find the Hasse--Witt matrix for this example as well.
Let $k$ be an algebraic closure of $\ff_{125}$.
By the strong approximation theorem, the vector space $H^1(X,\calO_X) \iso \calR/(\calR(0)+k(X))$
can be represented by the classes of r\'epartitions supported only at the point at infinity $\infty$
on the curve $X$. In fact, the r\'epartitions $r = \{r_P\}_{P\in X(k)}$ and $s = \{s_P\}_{P\in X(k)}$
defined by 
\[
r_P = \begin{cases} 2y/x & \text{if $P = \infty$;}\\ 0 & \text{otherwise}\end{cases}
\qquad\text{and}\qquad
S_P = \begin{cases} 2y/x^2 & \text{if $P = \infty$;}\\ 0 & \text{otherwise}\end{cases}
\]
give a basis for $\calR/(\calR(0)+k(X))$ that is dual to the basis $\{\omega_1,\omega_2\}$ of $H^0(X,\Omega_X^1)$ given in~\eqref{eqbasis} under the pairing \eqref{eqpairing}; we see this as follows. Let $z = x^2/y$, so that
$z$ is a uniformizing parameter for $X$ at $\infty$. We compute that
\begin{alignat*}{2}
\omega_1 &= \phantom{x}\, dx/y &&= \big(3z^2 + O(z^4)\big) \,dz\\
\omega_2 &= x\, dx/y           &&= \big(3 + 3z^2 + O(z^4)\big) \,dz\\
r_\infty &= 2y/x               &&= 2z^{-3} + 3z^{-1} + O(z)\\
s_\infty &= 2y/x^2             &&= 2z^{-1}.
\end{alignat*}
It follows easily that $(\omega_1,r) = (\omega_2,s) = 1$ and $(\omega_1,s) = (\omega_2,r) = 0$.

We compute also that
\begin{align*}
r_\infty^5 &= (2x^5 + 4x^4 + \alpha^2 x^3 + \alpha^{69} x^2 + \alpha^{77}x + \alpha^{94})y 
             + \alpha^{41} r_\infty + \alpha^{105} s_\infty + O(z)\\
s_\infty^5 &= 2y + 2r_\infty + \alpha^{95}s_\infty + O(z),
\end{align*}
and it follows that the Hasse--Witt matrix for our curve $X$ is given by
\[
A = \begin{pmatrix}
\alpha^{41} & 2\\
\alpha^{105} & \alpha^{95}
\end{pmatrix}.
\]
As expected, we see that $A$ is the transpose of Yui's matrix $Y$, that $B = (A^\tau)^\intercal$,
and that $A A^\sigma = \left(\begin{smallmatrix}
0&0\\0&0
\end{smallmatrix}\right).$

\subsection{A generalization.}
Garcia and Tafazolian generalize Manin and Yui's computation,
and calculate a matrix~\cite[p.~212]{GarciaTafazolian2008}
such that left-multiplication by this matrix 
gives the effect of the $n$-th iterate of the Cartier operator
in terms of the basis~$\calB$; the $(i,j)$ entry of their matrix is the
$p^n$-th root of the 
coefficient of $x^{ip^n-j}$ in the polynomial $f(x)^{(p^n-1)/2}$. 
The penultimate displayed equation on page~212 of their paper shows this matrix acting 
on the right, but the formul\ae{} presented elsewhere in in their paper make it clear that 
it acts on the left.

\section{Hasse--Witt matrices through the ages}
\label{sec:Manin&Yui}

As noted in the introduction, Hasse and Witt~\cite{HasseWitt1936} 
showed that various properties of a curve $X$ can be read off from the
action of Frobenius on $H^1(X,\calO_X)$, the equivalence classes
of r\'epartitions of the curve, and they 
associated a matrix to this semilinear operator.
In the paper in which he defined his operator on differential
forms, Cartier~\cite{Cartier1957} already noted a connection to
the Hasse--Witt matrix of the curve; Serre~\cite[\S~10]{Serre1958}
explains this well. Over the years, different authors have
made this connection more and more computationally explicit. 
In this section, we focus on the work of Manin and of Yui, because 
their papers are the ones referred to most often when present-day 
authors write about computational aspects of the Cartier operator.

\subsection{The work of Manin}
Manin published three works relevant to our discussion here. We treat them
each in turn.

In the first of these works~\cite{Manin1961}
(available also in an English translation~\cite{Manin1961english}),
Manin develops explicit formul\ae{} for computing the 
action of $\frobenius$ on $H^1(X,\calO_X)$. 
On one hand, the definition of the matrix $A$ in the second 
displayed equation on page~153 of~\cite{Manin1961} assumes a \emph{right} 
action.\footnote{The second displayed equation on page~245 
        of the English translation.
}
This is further emphasized in the first displayed equation on
page~154.\footnote{The final displayed equation on page~245
        of the English translation.
}
On the other hand, the change of basis formula in the last displayed 
formula on page 153, and the formula for the $g$-fold iterate of
$\frobenius$ on page 154, are valid provided matrices act on 
the \emph{left}.\footnote{The third displayed formula on page~245,
        and the $g$-fold iterate formula on the top of page~246,
        of the English translation.
}

The main result of this work 
(\cite[Theorem~1, p.~155]{Manin1961},~\cite[Theorem~1, p.~247]{Manin1961english})
considers a curve $X$ over a field with $q = p^e$ elements, and relates 
the characteristic polynomial of the Frobenius endomorphism of $\Jac X$
to the characteristic polynomial of a matrix representing the linear, 
$e$-fold iterate $\calF^e$.  The theorem as stated is correct, but only if we take
$A$ to be the matrix representing the Frobenius endomorphism of $H^1(X,\calO_X)$
\emph{acting on the left}.  However, since the matrix $A$ as defined in the
text before the theorem is taken to act on the \emph{right}, the theorem
is incorrect if it is taken in the larger context of the paper.

In the second paper we would like to discuss,
Manin~\cite{Manin1962,Manin1962english}
reconsiders some of these operators.  He works with the Cartier operator 
$\cartier$, observes that it is $\tau$-linear, and that it acts on the space
$H^0(X,\Omega^1_X)$.  He explicitly calculates a basis for the space of 
differentials on a particular hyperelliptic curve and computes
the action of the Cartier operator in terms of this basis,
using the same techniques that we reproduce here in 
Section~\ref{subsec:formulaire}. No matrices are written down,
so there are no obvious sign errors in this paper. Note, however,
that in this paper Manin considers the Cartier operator on $H^0(X,\Omega^1_X)$,
while in the preceding paper he considered the Frobenius operator on $H^1(X,\calO_X)$.

In Section~IV.5.2 of his paper on formal groups~\cite{Manin1963,Manin1963english},
Manin computes an operator that he \emph{calls} the Hasse--Witt
matrix --- and thus, in theory, should represent the action 
of $\frobenius$ on $H^1(X,\calO_X)$ --- but which actually
\emph{represents} the action of $\cartier$ on $H^0(X,\Omega^1_X)$,
as in the paper discussed in the preceding paragraph.
The formula Manin uses for iterates of this operator implicitly
(and incorrectly) assumes that it is $\sigma$-linear. 
This leads to errors in Section~IV.5.2; there are several problems 
with the displayed group of equations that deduce conditions on the
formal group of a curve's Jacobian from conditions on the 
equation of the curve (\cite[p.~86]{Manin1963}, \cite[p.~79]{Manin1963english}).
It seems to us that this paper may be the original source of a 
recurrent conflation in the literature of ``Hasse--Witt'' and 
``Cartier--Manin'' matrices.

\subsection{The work of Yui}
\label{subsecyui}

Yui~\cite{Yui1978} analyzes hyperelliptic curves with affine model $y^2=f(x)$,
and computes the Cartier operator $\cartier$ on $H^0(X,\Omega^1_X)$. 
(We remind the reader that Yui refers to the object we call the
\emph{Cartier operator} as the \emph{modified Cartier operator},
and that she denotes it by~$\cartier'$.)
In Theorem~2.1~\cite[p.~382]{Yui1978} and Theorem~2.2~\cite[p.~384]{Yui1978},
the formula for iterates is appropriate for a $\sigma$-linear operator, 
but $\cartier$ is $\tau$-linear.  Moreover, Lemma~D~\cite[p.~386]{Yui1978}
exploits the semilinear adjointness \eqref{eq:dual} between $\cartier$ and $\frobenius$, 
but overlooks the transpose necessary for such matrix calculations.  
Because of sign errors like these, Theorem~2.2~\cite[p.~384]{Yui1978} and
Lemma~E~\cite[p.~387]{Yui1978} are incorrect; the curve we discussed
in the Prologue gives a counterexample to both.

Although several explicit examples are worked out in Yui's paper, none of 
them can detect these inconsistencies.  Indeed, 
in Example~3.3~\cite[p.~391]{Yui1978} both $\cartier$ and $\frobenius$
are diagonalized by the natural basis, which hides ambiguity between left- and right-multiplication.  Moreover, both this example and 
Example~5.4~\cite[p.~400]{Yui1978} are worked out for curves over~$\ff_p$,
in which case $\sigma$- and $\tau$-linear operators are simply linear.

Yui writes at the end of the paper's introduction that the article
stemmed from her working through
Manin's papers~\cite{Manin1962,Manin1962english,Manin1963,Manin1963english},
so some of the sign errors in Yui's paper are reflections of Manin's earlier 
ambiguities between left actions and right actions and between $\sigma$-linear
and $\tau$-linear operators.
This paper also encourages the unfortunate conflation of the
concepts of the Hasse--Witt matrix and the Cartier--Manin matrix that began 
with Manin; we have already noted Lemma D~\cite[p.~386]{Yui1978}, which says
that the two matrices are ``identified'' with one another.

\section{Subsequent developments}
\label{sec:subsequent}

Explicit computational methods are becoming increasingly useful in 
arithmetic geometry, and this utility is reflected in the large number
of citations of the articles of Manin and Yui that we discussed in the 
preceding section.
Indeed, by consulting 
\href{http://www.ams.org/mathscinet/index.html}{MathSciNet} and the 
\href{http://webofknowledge.com}{Web of Science}, we found 
\the\numexpr\totvalue{appbibcounter}-1\relax 
 \ works
that refer to Yui's paper~\cite{Yui1978} or Manin's paper on Hasse--Witt 
matrices~\cite{Manin1961,Manin1961english}, and by personal knowledge we found one more. These works are listed below in a separate section of our bibliography.

It is somewhat worrisome to see so many citations, 
because --- as we have noted above --- 
these papers of Manin and Yui contain sign errors that
invalidate some  of their results.
To determine whether these sign errors have
propagated to other papers, we went through the \total{appbibcounter} articles
we found to see how they applied the results 
of Manin and Yui.
Of course, we could not go through all of these articles with great care;
for the most part, we limited ourselves to looking at how they made use of the work
of Manin and Yui described above, and it is possible we missed some subtleties.

In the vast majority of these works, we did not find any obvious
errors stemming from the citation of the papers of Manin and Yui.
For example:
\begin{itemize}
\item Sometimes the papers of Manin and Yui were given as general
references (for the computation of Hasse--Witt matrices or for something else),
and no particular results from the papers were used.
\nociteA{
AnbarBeelen2017,
Baker2000,
BauerJacobsonEtAl2008,
Bost2013,
CaisEllenbergEtAl2013,
Cassels1966,
CastryckStrengEtAl2014,
CornelissenOortEtAl2005,
DittersHoving1988,
DittersHoving1989,
Elkin2011,
ElkinPries2013,
Estrada-Sarlabous1991,
FarnellPries2013,
Galbraith2012,
GrangerHessEtAl2007,
Guerzhoy2005,
HashimotoMurabayashi1995,
IzadiMurty2003, 
Katz1972,
Kneppers1985,
KodamaWashio1988,
KodamaWashio1990,   
KudoHarashita2017,
Lennon2011,
MatsuoChaoEtAl2002,
MatsuoChaoEtAl2003,
Mazur1972,
Miller1976a,
Miller1976b,
Nygaard1981,
Nygaard1983,
Pries2005,
Pries2009,
PriesStevenson2011,
SarkarSingh2017,
Silverman2009,
SohnKim2009,
Takeda2011,
TakedaYokogawa2002,
Takizawa2007}
\item In some cases, specific results from Manin or Yui \emph{were} quoted, but
either they were not applied, or they did not contain any sign errors, 
or the sign errors were silently corrected.
\nociteA{
BalletRitzenthalerEtAl2010,
Ballico1994,
BassaBeelen2010,
BauerTeskeEtAl2005,
BouwDiemEtAl2004,
BuiumVoloch1995,   
CardonaNart2007,
Cherdieu2004,
Ditters1989,
DolgachevLehavi2008,
GlassPries2005,
Goodson2017a,
Goodson2017b,
Hawkins1987,
IbukiyamaKatsuraEtAl1986,
Miller1972,
Pacheco1990,
Ruck1986,
StohrVoloch1987,
Sung2017,
Washio1983,
Yui1980, 
Zapponi2008}
\item In some cases, statements containing sign errors
(quoted from Manin or Yui or elsewhere, or derived independently)
\emph{were} applied to specific examples, but in these examples the sign errors in the 
general formul\ae{} did not lead to errors in the specific cases. 
Incorrect formul\ae{} might not lead to errors, for example,
\begin{itemize}
\item if the genus of the curve is $1$; 
\item or, more generally, if the Hasse--Witt matrix is diagonal, so that $A$ commutes with all of its Galois conjugates;
\item or if the base field is $\ff_{p}$, so that no iteration is necessary;
\item or if the base field is $\ff_{p^2}$, so that $A \cdot A^\sigma = A\cdot A^\tau$; 
\item or in a number of other situations.
\end{itemize}
\nociteA{
Adolphson1980,
Alvarez2014,
Asada1988,
Cassou-NoguesChinburgEtAl2015,
FiteSutherland2016,
FurukawaKawazoeEtAl2004,
GhoshWard2015,
HarveySutherland2014,
HarveySutherland2016,
Hasegawa2013,
KodamaWashio1987,
Madden1978,
Olson1976,
Onishi2011,
Onishi2011english, 
Sohn2014,
Sullivan1975,
Tafazolian2012,
Ulmer1990,
Valentini1995,
Yui1988}
\end{itemize}

But in eight of these papers, incorrect results \emph{were} used in ways that
we felt required further investigation. We look at these papers here.

\subsection{Combining a theorem of Manin with a formula of Yui}
The paper of Gaudry and Harley~\citeA{GaudryHarley2000}, as well as 
the papers of Bostan, Gaudry, and 
Schost~\citeA{BostanGaudryEtAl2004,BostanGaudryEtAl2007}, 
all quote a result of Manin 
(\cite[Theorem~1, p.~155]{Manin1961}, \cite[Theorem~1, p.~247]{Manin1961english}; see also Section~\ref{subseczeta})
that relates the mod-$p$ reduction of the Weil polynomial of a
curve over $\ff_{p^e}$ to the characteristic polynomial of a matrix
\[ 
H_\pi = H H^{(p)}\cdots H^{(p^{e-1})},
\]
where $H$ is the Hasse--Witt matrix for the curve. As we noted earlier, Manin's
theorem is only correct as written if we take our matrices to act on the left.
However, the
papers of Bostan, Gaudry, Harley, and Schost under discussion take $H$ to be 
the matrix computed by Yui~\cite[p.~381]{Yui1978}.  Yui does intend for this matrix
to act on the left, but it represents the Cartier operator on differentials,
not the Frobenius operator on r\'epartitions, so Yui's matrix must be transposed
to give the Hasse--Witt matrix.  In other words, the na\"{\i}ve combination of Yui's 
matrix with Manin's theorem gives incorrect results.

This can be seen very concretely. 
Consider the genus-$2$ curve $X$ over $\ff_{27}$ defined
by $y^2 = x^5 + a^2 x^2 + a x$, where $a^3 - a + 1 = 0.$ 
On one hand, the matrices $H$ and $H_\pi$ from the cited papers are 
\[ 
H = \begin{bmatrix} a^2 & a\\1 & 0\end{bmatrix}
\textup{\quad and\quad}
H_\pi = H H^{(3)} H^{(9)} 
      = \begin{bmatrix} a^{12} & a^{14}\\a^{15} & a^{15}\end{bmatrix},
\]
and the characteristic polynomial $\kappa(t)$ of $H_\pi$ is
$t^2 + t + 1$. On the other hand, the characteristic polynomial 
of Frobenius for $X$ is $\chi(t) = t^4 + 6t^3 + 52t^2 + 162t + 729$, and it is visibly 
\emph{not} the case that $\chi(t) \equiv (-1)^2 t^2 \kappa(t) \bmod 3,$ as the
cited theorems claim. 

However, we suspect that Bostan, Gaudry, and Schost must have
implemented the computation of $H_\pi$ with the matrices in the opposite order
(or they transposed $H$, or something similar), because the example they 
present~\citeA[\S 5]{BostanGaudryEtAl2004} satisfies the basic 
sanity check that several randomly-chosen points on the Jacobian are
annihilated by the integer they give as the order of the Jacobian.

Likewise, Gaudry and Harley present an example~\citeA[\S 7.2]{GaudryHarley2000}
of a computation over $\ff_{p^4}$ in which they explicitly mention the 
order of the Jacobian modulo~$p$ computed by Manin's result, and the numerical
value they get shows that their computation must have involved either transposing
$H$ or computing $H_\pi$ with the factors reversed.

\subsection{Supersingular genus-\texorpdfstring{$2$}{2} curves}
\label{ss:genus2}
We found three papers that use Yui's computation of the
iterated Cartier operator
to determine when a genus-$2$ curve is supersingular.

Elkin~\citeA[\S 9]{Elkin2006} gives a characterization of supersingular 
genus-$2$ curves that includes a sign error.
This incorrect characterization does not affect the main part of
his work (for example, Theorems~1.1,~1.6, and~1.7~\citeA[pp.~54--56]{Elkin2006}), but we have not 
checked to see whether it affects the validity of his examples~\citeA[\S 9]{Elkin2006}.

Howe~\citeA{Howe2008} uses Yui's Lemma E~\cite[p.~387]{Yui1978} in the proof
of his Theorem~2.1~\citeA[p.~51]{Howe2008}, which claims that all supersingular 
genus-$2$ curves over a field of characteristic $3$ can be put into a certain
standard form. The proof as written is invalid, because the criterion for 
supersingularity has a sign error; however, the proof can easily be repaired
by using the correct criterion, and one can check that the theorem as stated
is true.

Zarhin~\citeA{Zarhin2004} also studies supersingular genus-$2$ curves in 
characteristic~$3$. In the proof of his Lemma~6.1~\citeA[p.~629]{Zarhin2004} he
correctly characterizes when a genus-$2$ curve is supersingular
in terms of a matrix that specifies the action of the Cartier operator.
Unfortunately, in a later paper~\cite[\S 5, p.~213]{Zarhin2005}
he provides a ``correction'' to
this proof that replaces the correct characterization with an incorrect one.  
Fortunately, this did not require changing the statement of the result he was
proving; the statement of
his Lemma~6.1~\citeA[p.~629]{Zarhin2004} is correct.

\subsection{Genus-\texorpdfstring{$3$}{3} curves of \texorpdfstring{$p$}{p}-rank \texorpdfstring{$0$}{0}}
We found one paper, by Elkin and Pries~\citeA{ElkinPries2007}, that uses
Yui's results to compute the moduli space of hyperelliptic
genus-$3$ curves of $p$-rank~$0$ in characteristic~$3$ and characteristic~$5$.
The notation in their Lemma~2.2 \citeA[p.~246]{ElkinPries2007}
is ambiguous, but when they apply this lemma in the proofs of Lemmas~3.3 
and~3.6~\citeA[pp.~248 and~250]{ElkinPries2007} they multiply the matrices in the wrong order. This 
invalidates their calculations of the defining equations of the 
moduli spaces. Pries reports that Theorem~4.2~\citeA[p.~251]{ElkinPries2007} still holds.

\subsection{Supersingularity versus superspeciality}
Yui's 1986 paper~\citeA{Yui1986} cites her 1978 paper~\cite{Yui1978}, 
as well as a paper of Nygaard~\citeA{Nygaard1981}, in the course of the
proof of Theorem~2.5~\citeA[p.~113]{Yui1986}. In particular, Yui cites these papers to
show that a curve over $\fbar_p$ has supersingular Jacobian (that is, 
its Jacobian is isogenous to a power of a supersingular elliptic curve)
if and only if the Cartier operator on its differentials is zero.  
In fact, Nygaard shows that the vanishing of the Cartier operator is 
equivalent to the Jacobian being \emph{superspecial} (that is, 
\emph{isomorphic} to a power of a supersingular elliptic 
curve)~\citeA[Theorem~4.1, p.~388]{Nygaard1981}. Furthermore, Yui herself 
gives examples showing that while the vanishing of the Cartier operator 
implies that the curve is supersingular, the converse is not 
true~\cite[Example~5.4, p.~400]{Yui1978}. 
Thus, Theorem~2.5~\citeA[p.~113]{Yui1986} is incorrect.

\section{Conclusion}

As we noted, most of the \total{appbibcounter} papers that cite Manin~\cite{Manin1961,Manin1961english}
or Yui~\cite{Yui1978} do not seem to have inherited any errors in their main results.
However, it might be prudent for authors who have used results from these \total{appbibcounter} papers
to double check that the results they quoted are indeed free of sign errors.

We conclude by repeating our supplication from the introduction:
Please be careful with terminology, and make a clear distinction between
the Cartier operator on differentials (represented by the Cartier--Manin 
matrix) and the Frobenius operator on $H^1(X,\calO_X)$ (represented by
the Hasse--Witt matrix). We hope that if such care is taken, there will
be no need in the future for another paper like this one.

\bibliographystyle{habbrvdoi}
\bibliography{cartier}

\bigskip
For the reader's convenience, we gather together here a list of all of
the papers that we are aware of that cite Manin's 1961 paper~\cite{Manin1961,Manin1961english} or 
Yui's 1978 paper~\cite{Yui1978}. We omit Yui's paper~\cite{Yui1978} itself,
even though it cites Manin~\cite{Manin1961english}.

\def\oldbibitem{} \let\oldbibitem=\bibitem
\def\bibitem{\stepcounter{appbibcounter}\oldbibitem}

\bibliographystyleA{habbrvdoi}
\bibliographyA{cartier}	

\bigskip

\end{document}